\documentclass{article}
\usepackage{amsfonts}
\usepackage{latexcad}

\advance\textheight by -\headheight
\advance\textheight by -\headsep
\newtheorem{proposition}{Proposition}

\newtheorem{theorem}{Theorem}
\newtheorem{hypoth}{Conjecture}
\newtheorem{coroll}{Corollary}

\def\eb{\bar{\bf 1}}                              
\def\rank{\mathop{{\rm rank}}\nolimits}           
\def\KT{K}                                        
\def\R{{\mathbb R}}                               
\def\ind{\mathop{\rm ind}\nolimits}               
\def\sqa{\sqcap\!\!\!\!\sqcup}                    
\def\epr{\hfill$\sqa$\medskip\par}                
\def\liml {\mathop{\lim}  \limits}                
\def\maxl {\mathop{\max}  \limits}                
\def\Im {\mathop{{\rm Im}}\nolimits}              
\def\l{\ell}                                      
\def\aa{\alpha}                                   
\def\G{\Gamma}                                    
\def\la{\lambda}                                  
\def\si{\sigma}                                   
\def\suml {\mathop{\sum}   \limits}               
\def\LT {\widetilde L}                            
\def\LcT {{\widetilde L_c}}                       
\def\J{\bar{J}}                                   
\def\dots{\cdots}                                 

\title{On the spectra of nonsymmetric Laplacian matrices\thanks{\sloppy
This work was supported by the Russian Foundation for Basic
Research under grant No.~02-01-00614}}

\author{Rafig Agaev\footnote{E-mail: {\tt arpo@ipu.ru}}$\:$
    and Pavel Chebotarev\footnote{%
Corresponding author.
E-mail: {\tt chv@lpi.ru; pchv@rambler.ru}}\\
{\small Trapeznikov Institute of Control Sciences of the Russian Academy of Sciences}\\
{\small 65 Profsoyuznaya Street, Moscow 117997, Russia}}

\date{}

\begin{document}
\maketitle

{\abstract\begin{minipage}[c]{42.5em}
A Laplacian matrix, $L=(\l_{ij})\in\R^{n\times n}$, has nonpositive off-diagonal entries and zero row sums.
As a matrix associated with a weighted directed graph, it generalizes the Laplacian matrix of an ordinary
graph. A standardized Laplacian matrix is a Laplacian matrix with $-{1\over n}\le\l_{ij}\le0$ at $j\ne i.$ We
study the spectra of Laplacian matrices and relations between Laplacian matrices and stochastic matrices. We
prove that the standardized Laplacian matrices $\LT$ are semiconvergent. The multiplicities of $0$ and $1$ as
the eigenvalues of $\LT$ are equal to the in-forest dimension of the corresponding digraph and one less than
the in-forest dimension of the complementary digraph, respectively. We localize the spectra of the
standardized Laplacian matrices of order $n$ and study the asymptotic properties of the corresponding domain.
One corollary is that the maximum possible imaginary part of an eigenvalue of $\LT$ converges to
${1\over\pi}$ as $n\to\infty.$
\end{minipage}}
\bigskip

\section{Introduction}

We consider the matrices $L=(\l_{ij})\in\R^{n\times n}$ such that
\begin{eqnarray}
\label{0504d1}               \l_{ij}&\le& 0,\quad j\ne i,  \\
\label{0504d2}   \sum_{j=1}^n\l_{ij}&=  & 0,\quad i=1,\ldots,n.
\end{eqnarray}

If, in addition, $L$ is symmetric and $\l_{ij}\in\{0,-1\},\;j\ne i,$
then $L$ is the {\em Laplacian matrix of an undirected graph\/}
(see, e.g.,~\cite{AndersonMorley}): $\l_{ij}=-1$ iff the graph has an edge between vertex $i$ and
vertex $j$, $j\ne i$; $\l_{ii}=$ the degree of vertex~$i$.
For the classical results on Laplacian matrices, we refer to
\cite{Kelmans65E,Kelmans66E,KelmansChelnokov74,AndersonMorley,GroneMerris90,GroneMerris94,Merris94,Merris95,Mohar92,Chung94}.

In the same way, every matrix of the form of (\ref{0504d1}) and (\ref{0504d2})
can be associated with a weighted directed graph~$\G$ whose arc weights are strictly positive:
\begin{equation}
\label{0504case1}
\l_{ij}=\left\{
\begin{array}{ll}
\mbox{minus the weight of arc }(i,j), & \hbox{if $j\ne i$ and $\G$ has an $(i,j)$ arc;} \\
0,                                    & \hbox{if $j\ne i$ and $\G$ has no $(i,j)$ arc;} \\
-\suml_{k\ne i}\l_{ik},               & \hbox{if $j=i$} \\
\end{array}
\right.
\end{equation}
(cf.\ the Kirchhoff matrix~\cite{Tutte}).

According to the famous matrix-tree theorem for digraphs \cite{Tutte}, the $(i,j)$ cofactor
of such a matrix $L(\G)$ is equal to the total weight of spanning
trees converging to $i$ in~$\G$, where the weight of a subgraph
is the product of the weights of its arcs. Generalizations of this theorem
have been obtained in \cite{FiedlerSedlacek58} and \cite{WKChen76,Chaiken82}.
We call the matrices that satisfy (\ref{0504d1}) and (\ref{0504d2}) {\em Laplacian matrices}.

It is easily seen \cite[p.~258]{CheAga02} that every Laplacian matrix is a singular M-matrix (for a
detailed discussion of M-matrices we refer to~\cite{BermanPlemmons}). By Ger\v{s}gorin's theorem (see,
e.g.,~\cite{HoJo86}), the real part of each nonzero eigenvalue of $L$ is strictly positive.

Laplacian matrices are closely connected with stochastic matrices. The basic relation is that the Laplacian
matrices can be defined as $\aa(I-P)$ matrices, where $\aa>0$, $I$ is the identity matrix, and $P$ is a
stochastic matrix. That is why the Laplacian matrices frequently appear in the theory of Markov
chains.\footnote{Studying the spectra of nonsymmetric Laplacian matrices is also needed for solving control
problems that involve graphs~\cite{FaxMurray02}.} Dmitriev and Dynkin \cite{DmitrievDynkin1} note that the
matrices $-L,$ where $L$ are Laplacian matrices, are the matrices of transition probability densities of
continuous time Markov chains. They strengthen the aforementioned result on the nonnegativity of the real parts
of the Laplacian eigenvalues by proving that for every eigenvalue $\la,$
\begin{equation}
\label{0504twolines}
-\left({\pi\over 2}-{\pi\over n}\right)\le{\rm arg}\,\la\le {\pi\over 2}-{\pi\over n}.
\end{equation}

In this paper, we further characterize the spectra of the Laplacian
matrices. To achieve a better comparability of the results for various $n$, we consider {\em standardized Laplacian
matrices}. These are the Laplacian matrices $\LT=(\l_{ij})$ with
\begin{equation}
\label{0504small}
-{1\over n}\le\l_{ij}\le0,\quad j\ne i,\:i=1,\ldots,n.
\end{equation}
For such matrices,
\begin{equation}
\label{0504smalldi}
0\le\l_{ii}\le{1-{1\over n}},\quad i=1,\ldots,n.
\end{equation}

If the class ${\mathcal G}_b$ of weighted digraphs with positive arc weights not exceeding
$b>0$ is considered and $L(\G)$ is the Laplacian matrix of a
weighted digraph $\G$ on $n$ vertices, then the {\it standardized
Laplacian matrix associated with\/ $\G$ in this class\/} is, by definition, $\LT(\G)=(nb)^{-1}L(\G).$

It follows from \cite[Theorem~1]{AndersonMorley} (or from the earlier results (2.7) and (2.10) of
\cite{KelmansChelnokov74}) that, in the case of symmetric Laplacian matrices, the eigenvalues of $\LT$ belong
to $[0,1]$. In the general case, (\ref{0504twolines}), (\ref{0504smalldi}), and Ger\v{s}gorin's theorem
immediately provide
\begin{proposition}
\label{130903l0}
The spectrum of a standardized Laplacian matrix belongs to the meet of the
closed disk with radius $1-1/n$ centered at $1-1/n$ and the closed
smaller angle bounded by the two half-lines drawn from $0$ through
the points\/ $e^{-({\pi \over 2}-{\pi \over n})i}$ and\/ $e^{({\pi
\over 2}-{\pi \over n})i}$.
\end{proposition}

For $n=7,$ Proposition~\ref{130903l0} is illustrated by Figure~1.
\begin{figure}[htb]
\begin{center}
\input 5.lp

Figure~1.
\end{center}
\end{figure}

In Section 2, we show that the eigenvalues of $\LT$ are also the eigenvalues of a certain stochastic matrix
and that $\LT$ is semiconvergent; the multiplicities of $0$ and $1$ as eigenvalues of $\LT$ are related with
the in-forest dimensions of the corresponding digraph $\G$ and the complementary digraph~$\G_c$,
respectively. In the subsequent sections, we advance in localizing the spectra of the nonsymmetric Laplacian
matrices.

\section{Laplacian matrices and stochastic matrices}

Let
\begin{equation}
\label{v_y_291003a}
\KT=I-\J,
\end{equation}
where $\J\in\R^{n \times n}$ is the matrix with all entries~$1/n$.
Obviously, $\KT\in\R^{n\times n}$ is the standardized Laplacian matrix of
the complete digraph with all arc weights $b$ in the class ${\mathcal G}_b$.

Define the matrices
\begin{equation}
\label{v_y_291003c}
P=\LT+\J
\end{equation}
and
\begin{equation}
\label{v_y_291003b}
\LcT=\KT-\LT.
\end{equation}

$P$ is stochastic, since it is nonnegative
with all row sums~1. In the class ${\mathcal G}_b$, $\LcT$ is the standardized Laplacian matrix
of the {\it complementary weighted digraph\/} $\G_c$ where the
weight attached to arc $(i,j),$ $j\ne i,$ is $b-w_{ij},$ $w_{ij}$ being the weight of this arc in~$\G$.
If $w_{ij}=b$ then $\G_c$ has no $(i,j)$ arc and vice versa: if $\G$ has no $(i,j)$ arc, then the weight of
$(i,j)$ in $\G_c$ is~$b$.
It follows from~(\ref{v_y_291003a}), (\ref{v_y_291003c}), and
(\ref{v_y_291003b}) that
\begin{equation}
\label{v_y_291003d}
P=I-\LcT.
\end{equation}

Let $\si(A)$ be the spectrum of a square matrix~$A$.

\begin{theorem}
\label{071103t1}
Suppose that $\LT$ is a standardized Laplacian matrix, $P$ is the stochastic matrix defined
by~$(\ref{v_y_291003c}),$ and $\LcT$ is the matrix defined by~$(\ref{v_y_291003b})$.
Then for $\la\not\in\{0,1\},$ the following statements are equivalent$:$\\
\indent {\rm (a)} $\la\in\si(\LT),\,$\\
\indent {\rm (b)} $\la\in\si(P),\,$ and\\
\indent {\rm (c)} $(1-\la)\in\si(\LcT),\,$\\
and these eigenvalues have the same geometric multiplicity.
Furthermore, $v$ is an eigenvector of $\LT$ corresponding to $\la\not\in\{0,1\}$ if and only if the
vector $x$ such that
\begin{equation}
\label{010504eq4}
x=\left(I-{\J\over{1-\la}}\right)v
\end{equation}
is an eigenvector of $P$ corresponding to~$\la$ and of $\LcT$ corresponding to $1-\la$.
\end{theorem}

{\bf Proof.}
Let $v$ be an eigenvector of $\LT$ corresponding to an
eigenvalue~$\la\not\in\{0,1\}$. Since $\LT\J=0$ and $\J^2=\J$, we have
$$
Px=
\left(\LT+\J\right)\left(I-{\J\over {1-\la}}\right) v=
 \LT v+\left(\J-{\J\over{1-\la}}\right)v=
 \la v-{{\la\J v}\over{1-\la}}=
\la \left(I-{\J\over {1-\la}}\right)v=\la x.
$$

Conversely, let $x$ be an eigenvector of $P$ corresponding
to~$\la\not\in\{0,1\}$. Observe that (\ref{010504eq4}) is equivalent to
\begin{equation}
\label{010504eq5}
v=\left(I-\frac{\J}{\la}\right)x.
\end{equation}
whenever $\la\not\in\{0,1\}$. Indeed, since $\J^2=\J,$
\[
\left(I-{\J\over{1-\la}}\right)\left(I-{\J\over{\la}}\right)=
I -{\J\over{\la}}-{\J\over{1-\la}}+{\J\over {\la(1-\la)}}=I.
\]

Using $\LT \J=0,$ we get
$$
\LT v=\LT\left(I-{\J\over\la}\right)x=\LT x=(P-\J)x=\la x-\J x=
\la\left(I-{\J\over\la}\right)x=\la v.
$$
The assertions about $\LcT$ follow from~(\ref{v_y_291003d}). The
geometric multiplicities of the eigenvalues are the same, since
the transformation (\ref{010504eq4}) is nondegenerate: (\ref{010504eq5}) defines its inverse.
\epr

\begin{theorem}
\label{0504tc}
Let $f_{\LT}(\la),$ $f_{P}(\la),$ and
$f_{\LcT}(\la)$ be the characteristic polynomials of $\LT,$ $P,$
and $\LcT,$ respectively. Then for all $\la\not\in\{0,1\},$
\begin{eqnarray}
\label{fP}
  f_{P}   (\la) &=& {\la-1\over\la}f_{\LT}(\la),\\
\label{fLcT}
  f_{\LcT}(\la) &=& (-1)^{n-1}{\la\over1-\la}f_{\LT}(1-\la).
\end{eqnarray}
\end{theorem}

{\bf Proof.} Let $c^{(1)},\ldots,c^{(n)}$ be the columns of the
characteristic matrix $\la I-\LT.$ Let $\eb=({1\over
n},\ldots,{1\over n})^T.$ Then, by (\ref{v_y_291003c}),
\begin{eqnarray}
\nonumber
  f_{P}(\la) &=& \det\left[c^{(1)}-\eb,\:c^{(2)}-\eb,\ldots,c^{(n)}-\eb\right]\\
\nonumber
   &=& \det\left[c^{(1)},\:c^{(2)},\ldots,c^{(n)}\right]\\
   &-& \det\left[\eb,c^{(2)},\ldots,c^{(n)}\right]
-\det\left[c^{(1)},\eb,\ldots,c^{(n)}\right]-\dots
-\det\left[c^{(1)},c^{(2)},\ldots,\eb\right].
\label{0405eq8}
\end{eqnarray}

Note that
$\det\left[c^{(1)},c^{(2)},\ldots,c^{(n)}\right]=f_{\LT}(\la)$ and
$\eb=(\la n)^{-1}(c^{(1)}+\dots+c^{(n)}),$ provided that $\la\ne0$. Therefore
$$
\det\left[\eb,c^{(2)},\ldots,c^{(n)}\right]=
\det\left[{1\over\la n}\left(c^{(1)}+\sum_{k=2}^n c^{(k)}\right),c^{(2)},\ldots,c^{(n)}\right]
={1\over\la n}f_{\LT}(\la).
$$
This expression works for each of the last $n$ items of~(\ref{0405eq8}).
Substituting this in (\ref{0405eq8}) gives
$f_P(\la)\!=\!f_{\LT}(\la)-n{1\over\la n}f_{\LT}(\la)\!=\!{\la-1\over\la}f_{\LT}(\la).$ To prove (\ref{fLcT}),
note that by (\ref{v_y_291003d}), for all $\la\ne1,$
$f_{\LcT}(\la)\!=\!\det(\la I-I+P)\!=\!(-1)^n\det(I-\la I-P)\!=\!(-1)^n f_P(1-\la)\!
=\!(-1)^{n-1}{\la\over1-\la}f_{\LT}(1-\la)$.
\epr

{\bf Remark 1.} Since (\ref{fP}) and (\ref{fLcT}) are true for all $\la$, except for
$0$ and $1$, they become valid for all~$\la$ after reducing the common factors $\la$ and $1-\la$, respectively.
\smallskip

An analogy of (\ref{fLcT}) for symmetric matrices is due to Kelmans \cite[Property~5]{Kelmans65E}.

Recall that a square matrix $A$ is {\em semiconvergent\/} if $\lim_{k\to\infty}A^k$ exists. The {\it index\/}
of $A$, $\ind A$, is the least $k=0,1,\ldots,n-1$ such that $\rank A^{k+1}=\rank A^k.$

\begin{theorem}
\label{100903t1}
For every standardized Laplacian matrix $\LT$ and the corresponding stochastic matrix $P$ defined
by $(\ref{v_y_291003c}),$ $\LT$ and $P$ are semiconvergent.
\end{theorem}

{\bf Proof.} Let $\rho(A)$ be the spectral radius of a matrix~$A$.
A square matrix $A$ is semiconvergent (see, e.g.,
\cite[Problem~4.9]{BermanPlemmons}) if and only if it satisfies:

(a) $\rho(A)\le 1$,

(b) If $\la\in\si(A)$ and $|\la|=1$, then $\la=1$, and

(c) $\rank(I-A)^2=\rank(I-A)$, i.e., $\ind(I-A)\in\{0,1\}$.

We now prove that $P$ satisfies (a), (b), and (c).
(a)~results from the stochasticity of~$P$.

(b)~Since all diagonal entries of the singular M-matrix $\LcT=I-P$
are less than $1$, $|\la|=1$ implies $\la=1$ (see, e.g., \cite[page
153]{BermanPlemmons}).

(c)~Since the index
of any Laplacian matrix is~1 \cite[Proposition~12]{CheAga02},
we have $\rank(I-P)^2=\rank(I-P)$, and $P$
satisfies~(c).
Thus, $P$ is semiconvergent.

Finally, it follows from (\ref{v_y_291003c}) and $\LT\J=0$ that
$\LT=P-\J,$
$\LT^2=\LT(P-\J)=\LT P=(P-\J)P,$
$\ldots,$
$\LT^k=\LT(P-\J)P^{k-2}=\LT P^{k-1}=(P-\J)P^{k-1},\ldots.$

Therefore,
\begin{equation}
\label{limus}
\lim_{k\to\infty}\LT^k=(P-\J)\lim_{k\to\infty}P^{k-1},
\end{equation}
and $\LT$ is semiconvergent as well.
\epr

Let $m_A(\la)$ be the algebraic multiplicity of $\la$ as an
eigenvalue of~$A$; let $V_A(\la)$ be the set of eigenvectors of $A$ corresponding to~$\la$.
An eigenvalue is called {\em semisimple\/} if its algebraic and geometric multiplicities coincide.
The {\it in-forest dimension\/} of a digraph is the minimum possible number
of trees in a spanning converging forest (also called {\it in-forest}) of the digraph
\cite[p.~255]{CheAga02}.

\begin{theorem}
\label{100903t5}
Let $d$ and $d_c$ be the in-forest dimensions of the digraph
corresponding to~$\LT$ and the complementary digraph,
respectively. Then:\\
{\rm(i)} $\:\:\:\:\:m_{\LT}(0)=d,$ $\quad\quad\: m_{\LT}(1)=d_c-1;$\\
{\rm(ii)} $\:\:\:m_{P}(0)=d-1,$ $\:\:m_{P}(1)=d_c;$\\
{\rm(iii)} $m_{\LcT}(1)=d-1,$ $\:m_{\LcT}(0)=d_c,$\\
and these eigenvalues are semisimple.\\
{\rm(iv)} If $v\in V_{\LT}(0)$ and $\KT v\ne0,$ then $\KT v\in V_P(0)=V_{\LcT}(1);$\\
\phantom{.}\quad$\:\:$if $x\in V_P(1)=V_{\LcT}(0)$ and $\KT x\ne0,$ then $\KT x\in V_{\LT}(1)$.
\end{theorem}

{\bf Proof.} By \cite[Proposition~12]{CheAga02}, $\rank L=n-d$ and $\ind L=1,$ therefore $m_{\LT}(0)=d$ and
$0$ is a semisimple eigenvalue. By Theorem~\ref{100903t1}, $\LT$ is semiconvergent. Thus, $\LT$ satisfies
condition (c) in the proof of Theorem~\ref{100903t1}, which implies the semisimplicity of~1, provided that
$1\in\si(\LT)$. Indeed, by (c), the maximum order of the Jordan blocks corresponding to~1 in the Jordan form
of $\LT$ is at most~1. By virtue of (\ref{fLcT}) and Remark~1, $m_{\LT}(1)=m_{\LcT}(0)-1$. The statements
regarding $\LcT$ follow from those about $\LT;$ after that, (ii) and the assertions regarding $P$ result
from~(\ref{v_y_291003d}). The proof of (iv) is straightforward.
\epr

For the case of symmetric Laplacian matrices, the identities $m_{\LT}(0)=d$ and $m_{\LT}(1)=d_c-1$ are due to
Kelmans \cite{Kelmans65E}. In that case, $d$ and $d_c$ reduce to the numbers of components of the
corresponding graphs.

In the following sections, we focus on the eigenvalues with nonzero imaginary parts.

\section{A region that contains the Laplacian spectra}

\begin{theorem}
\label{130903t2}
Let $\LT$ be a standardized Laplacian matrix of order~$n$. Then the spectrum of $\LT$ belongs to the meet of:

-- two closed disks, one centered at $1/n$, the other centered at
$1-1/n$, each having radius $1-1/n$,

-- two closed smaller angles, one bounded with the two half-lines
drawn from $1$ through $e^{-2\pi i/n}$ and $e^{2\pi i/n}$, the other
bounded with the half-lines drawn from $0$ through
$e^{-({\pi\over 2}-{\pi \over n})i}$ and
$e^{({\pi \over 2}-{\pi \over n})i}$, and

-- the band $|\Im(z)|\le{1\over 2n}\cot{\pi\over 2n}.$
\end{theorem}

{\bf Proof.} In addition to Proposition~\ref{130903l0}, note that the domain specified in that proposition
contains the spectrum of $\LcT$ too. Therefore, the spectrum of $P=I-\LcT$ and, by Theorem~\ref{071103t1},
the spectrum of $\LT$ belongs to the disk centered at $1/n$ and the angle with vertex at $1$, as specified in
Theorem~\ref{130903t2}.

Furthermore, according to Pick's theorem \cite{Pick1922}, if $\la$ is an eigenvalue of $A\in\R^{n \times n}$,
then $|\Im(\la)|\leq g\,\cot(\pi/2n)$, where $g=\maxl_{k,j}|c_{kj}|$ and $C=(c_{kj})={1\over 2}(A-A^*)$. If
$A=\LT,$ then $g\le{1\over 2n}$ by virtue of~(\ref{0504small}). This completes the proof.
\epr

{\bf Remark 2.} It follows from the stochasticity of $I-\LT$ and \cite[Theorem~4]{DmitrievDynkin2} that
$\si(\LT)$ contains an eigenvalue with argument ${\pi\over2}-{\pi\over n}$ if and only if the corresponding
digraph is a Hamiltonian cycle on $n$ vertices. In this case, such an eigenvalue $\la$ is unique,
$|\la|\le{2\over n}\sin{\pi\over n}$, and $\Im(\la)\le{1\over n}\sin{2\pi\over n}$. These bounds are reached
when $\LT={1\over n}(I-Q),$ where $Q$ is the standard circulant permutation matrix (see the following section).
The components of any
corresponding eigenvector are the vertices of a regular polygon. Similarly, $\si(\LT)$ contains an eigenvalue
that belongs to $[1,e^{{2\pi i\over n}}]$ if and only if the complementary digraph $\G_c$ is a Hamiltonian
cycle on $n$ vertices. As above, such an eigenvalue $\la'$ is unique and $\Im(\la')\le{1\over
n}\sin{2\pi\over n}$.

Theorem~\ref{130903t2} and Remark~2 are illustrated by Figure~2, where $n=7$.

\addvspace{25mm}
\begin{center}
\input 41.lp
\phantom{....................}

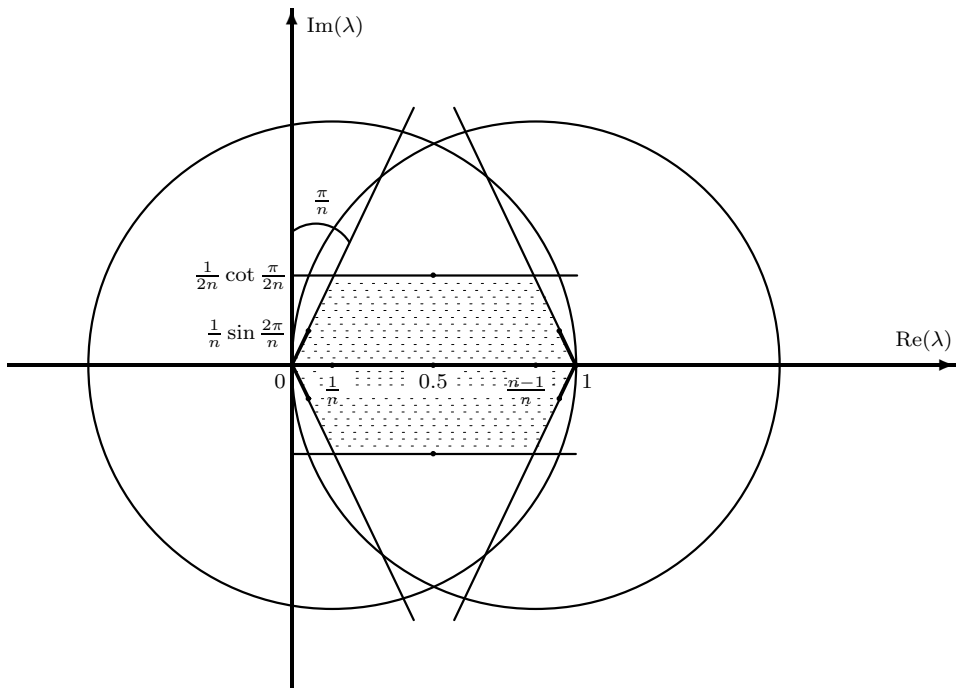
Figure~2.
\end{center}

The region $R$ specified by Theorem~\ref{130903t2} (hatched in Figure~2) is a hexagon whenever $4\le n\le
18,$ however its boundary contains arcs when $n>18$. This region reduces to the $[0,1]$ segment when $n=2$
and to the rhombus
with vertices $0,(0.5,(2\sqrt{3})^{-1}),1,$ and $(0.5,-(2\sqrt{3})^{-1})$
when $n=3.$

\section{The polygon of Laplacian eigenvalues}

In this section, we introduce a polygon that consists of Laplacian eigenvalues. As above, the order $n$ of
standardized Laplacian matrices is fixed.

Let $Q=(q_{kj})$ be the standard $n\!\times\!n$ circulant permutation matrix:
$q_{kj}=1$ if $j-k\in\{1,1-n\}$ and $q_{kj}=0$ otherwise. Consider $n-1$ matrices
\begin{equation}
\label{160704eq1}
L_k={1\over n}\left(k I-Q-Q^2-\dots-Q^k\right),\quad k=1,\ldots,n-1.
\end{equation}

Since $Q^i,\:i=1,\ldots,n-1,$ are permutation matrices whose patterns of ones are different entrywise, the
$L_k$'s have zero row sums, and their off-diagonal entries belong to $[-\frac{1}{n},0]$. Thus, they are
standardized Laplacian matrices.

The characteristic polynomial of $Q$ is $f_Q(\mu)=\mu^n-1,$ therefore the spectrum of $Q$ contains
\begin{equation}
\label{190704eq1}
\mu=e^{-2\pi i/n}.
\end{equation}
Consequently, the spectrum of $L_k$ contains
\begin{eqnarray}
\label{270604eq2}
\la_k&=&{1\over n}(k-\mu-\mu^2-\dots-\mu^k)\\
\label{270604eq2a}
      &=&{1\over n}\left(k-e^{-2\pi i/n}-\dots-e^{-2k\pi i/n}\right),
\end{eqnarray}
as well as $\bar\la_k$.

By means of the identity
$$
e^{\aa i}+e^{2\aa i}+\dots+e^{k\aa i}
={\sin{k\aa\over 2}\over\sin{\aa\over 2}}\,e^{(k+1)\aa i/2},
$$
(\ref{270604eq2a}) simplifies to:
\begin{equation}
\label{270604eq3}
\la_k
=n^{-1}\left(k
      -{\sin{k\pi\over n}\over
        \sin{\pi  \over n}}\,
      e^{-(k+1)\pi i/n}
       \right)
={k\over n}
      -{\sin{k\pi\over n}\over
       n\sin{\pi  \over n}}
  \left(\cos{(k+1)\pi\over n}-
       i\sin{(k+1)\pi\over n}\right).
\end{equation}
Let $S$ be the convex polygon with $2(n-1)$ vertices:
$\la_0=0,\,\la_1,\ldots,\la_{n-2},\,\la_{n-1}=1,\,\bar{\la}_{n-2},\ldots,\bar{\la}_1$.

\begin{theorem}
\label{270604t1}
Every point of the polygon $S$ is an eigenvalue of some standardized Laplacian matrix of order~$n$.
\end{theorem}

{\bf Proof.} As shown above, every vertex $\la_k$ of $S$ is the eigenvalue of the corresponding~$L_k$. Every
other point $s\in S$ (in what follows, we suppose $\Im(s)\ge0$, since the spectrum of a real matrix is
symmetric with respect to the $x$-axis) can be represented by a convex combination $f(0,\la_k,\la_{k+1})$ of
three vertices of $S$: $0$ and the two neighboring vertices $\la_k$ and $\la_{k+1}$ such that the half-line
$[0,s)$ meets the segment $[\la_k,\la_{k+1}].$ By definition, $\la_k$ and $\la_{k+1}$ are the values of the
polynomials (\ref{270604eq2}) for $\mu$ (we denote these polynomials by $p(\cdot)$ and $q(\cdot)$,
respectively); furthermore, they are the eigenvalues of the matrices $L_k=p(Q)$ and $L_{k+1}=q(Q)$,
respectively. Thus, $s=f(0,p(\mu),q(\mu))$, where $f(0,p(\cdot),q(\cdot))$ is a polynomial. Consequently, $s$
is the eigenvalue of the matrix $f(0,p(Q),q(Q)),$ which is a standardized Laplacian matrix as a convex
combination of such matrices.
\epr

It is easily seen that the polygon $S$ reduces to $[0,1]$ when $n=2$. For $n=3,$ $S$ is the rhombus with
vertices $0,\,(0.5,(2\sqrt{3})^{-1}),\,1,\,(0.5,-(2\sqrt{3})^{-1}),$ and it coincides with the region
described in Theorem~\ref{130903t2} as well.

Let
\begin{equation}
\label{010404eq1}
z(n)=\sup\{|\Im(\la)|:\la\mbox{\ is an eigenvalue of a standardized\ }n\!\times\!n\mbox{\ Laplacian matrix}\}.
\end{equation}

By Theorem~\ref{130903t2}, for every $n=2,3,\ldots,\;$ $z(n)\le{1\over 2n}\cot{\pi\over 2n}.$

\begin{proposition}
\label{210404p1}
If $n$ is odd, then
\begin{equation}
\label{280404eq1}
z(n)={1\over 2n}\cot{\pi\over 2n},
\end{equation}
moreover, $z(n)=\Im(\la_{(n-1)/2}),$ where $\la_{(n-1)/2}$ is defined by~$(\ref{270604eq2})$.
\end{proposition}

The proposition is proved by substituting $k=(n-1)/2$ in~(\ref{270604eq3}).

If $n>2$ is even, then the upper bound ${1\over 2n}\cot{\pi\over 2n}$ is not reached by the imaginary parts
of the vertices of $S$; more specifically, the following proposition is easy to verify.
\begin{proposition}
\label{210404p1a}
If $n>2$ is even, then
$$
\max_{0\le k\le 2n-1}\Im(\la_k)=\Im(\la_{n/2})={1\over n}\cot{\pi\over n}<{1\over 2n}\cot{\pi\over 2n}.
$$
\end{proposition}

The following conjecture is yet unproved.

\begin{hypoth}
The set of all eigenvalues of standardized Laplacian matrices of order $n$ coincides with the polygon~$S$.
\end{hypoth}

\section{Some asymptotic properties of the Laplacian spectra}

\begin{theorem}
\label{200704t1}
The boundary of the polygon $S$ converges$,$ as $n\to\infty,$ to the curve made up by the
parts of two cycloids whose parametric equations are:
$\:z(\tau)=x(\tau)+i\,y(\tau)\:$ and $\:z(\tau)=x(\tau)-i\,y(\tau),$ where
\begin{eqnarray}
\label{190704eq2}
   x(\tau)&=&(2\pi)^{-1}(\tau-\sin\tau),\\
\label{190704eq3}
   y(\tau)&=&(2\pi)^{-1}(1-\cos\tau),
\end{eqnarray}
and $\tau\in[0,2\pi]$.
\end{theorem}

{\bf Proof.} For each $\tau\in[0,2\pi),$ consider any positive integer sequence $k(n)$ such that
$\liml_{n\to\infty}{k(n)\over n}={\tau\over 2\pi}.$ Then, by virtue of (\ref{270604eq3}),
\[
\liml_{n\to\infty}\la_{k(n)} ={1\over 2\pi}\left(\tau-2\sin{\tau\over 2}\,\cos{\tau\over
2}+2i\sin^2{\tau\over 2}\right) ={1\over 2\pi}(\tau-\sin\tau+i(1-\cos\tau))
.
\]
Taking the real and imaginary parts of this expression, we obtain the desired statement.
\epr

Figure~3 
shows the polygons $S$ at $n=4$ and $n=5,$ along with the limit curve specified in Theorem~\ref{200704t1}.

\begin{figure}[htb]
\begin{center}
\label{240704f3}
\input cycloid.lp

Figure~3.
\end{center}
\end{figure}

\begin{coroll}
$\liml_{n\to\infty}z(n)=1/\pi$, and the convergence is from below.
\end{coroll}

This corollary is easily proved using Propositions~\ref{210404p1} and~\ref{210404p1a}.
\bigskip

\noindent{\bf Acknowledgments.} We are grateful to two anonymous referees for their helpful comments.

\end{document}